\documentclass[english]{elsarticle}
\usepackage[T1]{fontenc}
\usepackage[latin9]{inputenc}
\usepackage{float}
\usepackage{amsmath}
\usepackage{amsthm}
\usepackage{amssymb}
\usepackage{graphicx}

\makeatletter
\usepackage{url}

\makeatother

\usepackage{babel}
\begin{document}
\title{A comparison of semi-Lagrangian discontinuous Galerkin and spline based Vlasov solvers in four dimensions\tnoteref{label1}}
\tnotetext[label1]{The computational results presented have been achieved [in part] using the Vienna Scientific Cluster (VSC).}
\author[tub,uibk]{Lukas Einkemmer\corref{cor1}} \ead{einkemmer@na.uni-tuebingen.de}
\address[tub]{Department of Mathematics, University of T\"ubingen, Germany}
\address[uibk]{Department of Mathematics, University of Innsbruck, Austria}
\cortext[cor1]{Corresponding author}
\begin{abstract}
The purpose of the present paper is to compare two semi-Lagrangian methods in the context of the four-dimensional Vlasov--Poisson equation. More specifically, our goal is to compare the performance of the more recently developed semi-Lagrangian discontinuous Galerkin scheme with the de facto standard in Eulerian Vlasov simulation (i.e.~using cubic spline interpolation). To that end, we perform simulations for nonlinear Landau damping and a two-stream instability and provide benchmarks for the SeLaLib and sldg codes (both on a workstation and using MPI on a cluster).

We find that the semi-Lagrangian discontinuous Galerkin scheme shows a moderate improvement in run time for nonlinear Landau damping and a substantial improvement for the two-stream instability. It should be emphasized that these results are markedly different from results obtained in the asymptotic regime (which favor spline interpolation). Thus, we conclude that the traditional approach of evaluating numerical methods is misleading, even for short time simulations.  In addition, the absence of any All-to-All communication in the semi-Lagrangian discontinuous Galerkin method gives it a decisive advantage for scaling to more than 256 cores.
\end{abstract}  
\begin{keyword} semi-Lagrangian discontinuous Galerkin method, Vlasov--Poisson equation, semi-Lagrangian methods, performance comparison, parallelization \end{keyword}
\maketitle

\section{Introduction}

Many problems in computational plasma physics require the solution
of kinetic equations. Due to the high dimensionality of these models,
particle methods have been, and still are, heavily used to conduct
these simulations. However, due to the rapid increase in computational
performance, numerical methods based on the discretization of the
entire phase space (the so-called Eulerian approach) have gained increasing
attention in recent years (at least for four-dimensional and certain
five-dimensional problems). These Eulerian methods offer superior
accuracy and the possibility to resolve regions of the phase space
with low density, but are extremely demanding from a computational
point of view. Consequently, it is important to develop efficient
numerical methods that can be parallelized to state of the art supercomputers.

While standard techniques from fluid dynamics can be applied to kinetic
equations, these usually suffer from a severe restriction in the time
step size due to a Courant\textendash Friedrichs\textendash Lewy (CFL)
condition. On the other hand, the use of implicit methods would significantly
increase the memory footprint, which is already a scarce resource
for high-dimensional simulations. This is an even more pronounced
issue for accelerators such as graphic processing units (GPUs) and
the Intel Xeon Phi (which in all likelihood will become more pervasive
in the future). Semi-Lagrangian methods offer a viable alternative.
These methods follow the characteristics backward in time. In many
kinetic models (such as the Vlasov\textendash Poisson and Vlasov\textendash Maxwell
equations) applying a splitting procedure results in sub-steps for
which the characteristics can be determined analytically (if this
is not the case an explicit ordinary differential equation solver
is usually used). For the Vlasov\textendash Poisson equation this
approach has been first suggested in the seminal paper of Cheng \&
Knorr \cite{cheng1976}. Later it was extended in various ways to
the Vlasov\textendash Maxwell equation \cite{crouseilles2015hamiltonian,valis,mangeney2002}.
However, since the endpoint of a characteristic curve does not necessarily
coincide with the grid used, an interpolation procedure has to be
employed. An attractive choice is to reconstruct the desired function
by spline interpolation (see, for example, \cite{sonnendrucker1999,filbet2003}),
which according to \cite{sonnendruecker2011} is still considered
the de facto standard in Vlasov simulations. The main downside of
this approach is that a tridiagonal linear system of equations has
to be solved to construct the spline. While this can be done very
efficiently in a sequential implementation, for modern supercomputers
(which require scaling to hundreds of thousands of cores) such global
dependencies seriously limit the utility of this algorithm.

On the other hand, the semi-Lagrangian discontinuous Galerkin method
employs a piecewise polynomial approximation in each cell of the computational
domain (see, for example, \cite{rossmanith2011,crouseilles2011,einkemmer2014,qiu2011}).
In the case of an advection equation the discretized function is translated
and then projected back to the appropriate subspace of piecewise polynomial
functions. This method, per construction, is mass conservative and
only accesses two adjacent cells in order to compute the necessary
projection (this is true independent of the order of the approximation).
The local nature of this method also makes it easy to incorporate
block structured mesh refinement and employing limiters to maintain
positivity. Positivity limiters were also investigated for spline
interpolation (see \cite{zerroukat2005,zerroukat2006,crouseilles2010}),
but these schemes are much more involved. For the semi-Lagrangian
discontinuous Galerkin scheme mathematical rigorous convergence results
are available in \cite{einkemmer2014convergence,einkemmer2014}. The
method is fully explicit (i.e.~no linear system has to be solved
to advance the solution in time) and thus it is easier to implement
(especially on parallel architectures) and shows a more favorable
communication pattern. We should note that some measures have been
taken to improve the parallel scalability of cubic spline interpolation
\cite{crouseilles2009}. However, even for this approach a relatively
large communication overhead is incurred. This is due to the fact
that the boundary condition for the local spline reconstruction requires
a large stencil if the desirable properties of the global cubic spline
interpolation are to be preserved (the method derived in \cite{crouseilles2009}
requires a centered stencil of size $21$). In addition, a very efficient
mixed precision implementation is available for the semi-Lagrangian
discontinuous Galerkin method \cite{einkemmer2016}. This is particular
useful for modern accelerators (such as GPUs) and helps to alleviate
the problem of memory scarcity on such systems. Furthermore, it is
worth mentioning that a filamentation filtration strategy has been
developed \cite{einkemmer2014recurrence} (for problems where numerical
recurrence is an issue) and that the error propagation is superior
to other semi-Lagrangian schemes \cite{steiner2013semi,einkemmer2015error}
(although to our knowledge no proof of this fact is available for
the Vlasov equation).

Before proceeding, let us note that a range of other methods have
been used to solve the Vlasov equation as well. A method that is related
to the discontinuous Galerkin approach and which has been widely employed
in plasma simulations is the van Leer scheme (see, for example, \cite{fijalkow1999,mangeney2002,galeotti2005,califano2006}).
What distinguishes the discontinuous Galerkin method from the van
Leer scheme is that the coefficients in the corresponding basis expansion
are stored directly in computer memory. In contrast, the van Leer
scheme replaces these coefficients by performing an approximation
using suitable differences on an equidistant grid. Its implementation
and performance characteristics are thus more closely related to other
finite difference schemes. 

Despite the many advantages of the semi-Lagrangian discontinuous Galerkin
method outlined above, there is certainly no universal agreement that
this method will be a viable alternative for large scale high-dimensional
kinetic simulation. Specifically, it is often argued, both by reviewers
and in informal conversation, that these methods require more degrees
of freedom to achieve the same precision (compared to say spline interpolation)
and that this is not maintainable in a higher-dimensional setting.
This argument is certainly not easily dismissed. Let us assume that
instead of $n$ degrees of freedom per direction our numerical approach
requires $\alpha n$, with $\alpha>1$. Then the memory footprint
and run time are increased by a factor $\alpha^{d}$, where $d$ is
the dimensionality of the problem. For example, even a moderate increase,
for example $\alpha=1.2$, results in a substantial increase of memory
footprint/run time (a factor of $2$ for $d=4$ and factor of $3$
for $d=6$).

Of course, the important question then becomes if the semi-Lagrangian
discontinuous Galerkin method actually requires more degrees of freedom
and if so by how much. Although this is a common problem for classic
(i.e.~non semi-Lagrangian) discontinuous Galerkin methods and has,
for example, lead to so-called hybrid schemes (see, for example, \cite{lehrenfeld2016high}),
the situation considered here is quite different. Specifically, the
problem is reduced to a sequence of one-dimensional advection equations
which are then solved using the method of characteristics. In the
present paper we will investigate if the semi-Lagrangian discontinuous
Galerkin method requires more degrees of freedom compared to spline
interpolation. This will be done in the context of the Vlasov\textendash Poisson
equation \begin{equation}\label{eq:vlasov-poisson}
\begin{aligned}
& \partial_{t}f(t,x,v)+v\cdot\nabla_{x}f(t,x,v)+E(f)(x)\cdot\nabla_{v}f(t,x,v)=0 \\
& \nabla\cdot E(f)(x)=\int f(t,x,v)\,\mathrm{d}v-1,\qquad\;\;\nabla\times E(f)(x)=0.
\end{aligned}
\end{equation}However, before we proceed and discuss the comparative efficiency
of the semi-Lagrangian discontinuous Galerkin method (section \ref{sec:performance})
and the four-dimensional simulations conducted (sections \ref{sec:nl}
and \ref{sec:tsi}), let us dwell a bit more on the argument (against
the discontinuous Galerkin method) outlined above. To do that we will
consider nonlinear Landau damping as a test problem and, for the time
being, still restrict ourselves to the two-dimensional case (for reasons
that will be apparent shortly; the four-dimensional variant of this
problem is discussed in section \ref{sec:nl}). That is, we consider
the initial value
\[
f(0,x,v)=\frac{1}{\sqrt{2\pi}}(1+\epsilon\cos(kx))\mathrm{e}^{-\frac{1}{2}v^{2}},
\]
where $\epsilon=\tfrac{1}{2}$ and $k=\tfrac{1}{2}$, on the domain
$[0,4\pi]\times[-6,6]$. A splitting procedure is employed to obtain
a sequence of advection equations which are discretized using either
a semi-Lagrangian approach based on spline interpolation or the semi-Lagrangian
discontinuous Galerkin approach (both numerical methods will be discussed
in some detail in section \ref{sec:numerical-method}). The error
committed by these two numerical schemes is computed by comparing
it to a reference solution (with a sufficiently fine resolution).

\begin{figure}
\begin{centering}
\includegraphics[width=12cm]{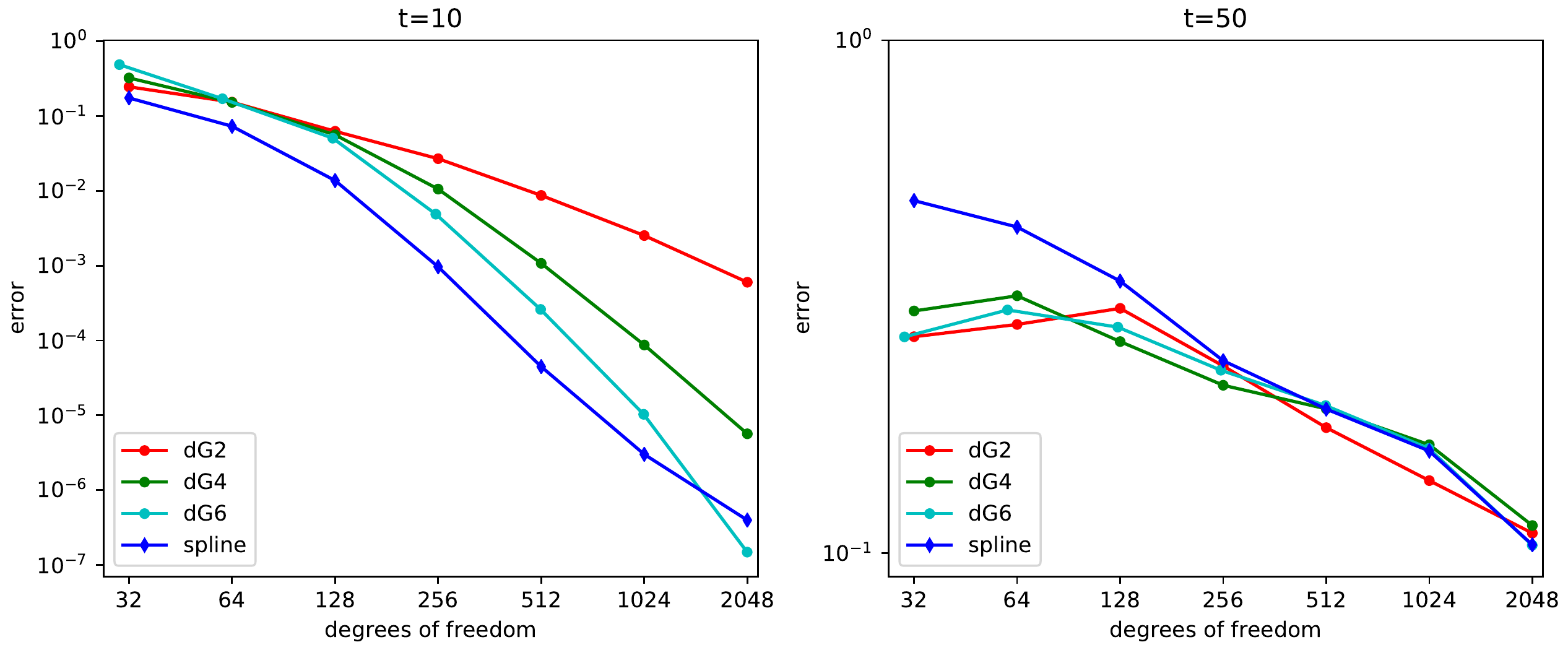}
\par\end{centering}
\caption{The error, as a function of the number of degrees of freedom (per
direction), for nonlinear Landau damping is shown at $t=10$ (left)
and $t=50$ (right). The error is computed in the infinity norm and
is normalized by the norm of the reference solution. The reference
solution is computed using spline interpolation with $8192$ grid
points (per direction). For all simulations a Strang splitting procedure
with time step size $\tau=0.1$ has been used. \label{fig:order}}
\end{figure}

Many publications (see, for example, \cite{einkemmer2017study,crouseilles2011,cheng2014discontinuous}),
including work by the authors of this paper, have reported such results.
To facilitate the present discussion. we have computed the error in
the particle-density $f$ at time $t=10$ as a function of the degrees
of freedom (per direction). These results are shown on the left in
Figure \ref{fig:order}. We clearly observe that, at best, the discontinuous
Galerkin method struggles in a direct comparison with spline interpolation.
The sixth-order discontinuous Galerkin method (dG6) requires approximately
50\% more degrees of freedom to obtain the same accuracy. Thus, it
seems that the discussion has come to a close and the criticism outlined
above seems to be completely justified. However, the major caveat
here is that in practical applications no one is interested in solving
nonlinear Landau damping only up to $t=10$. Now, if we look at the
error at $t=50$ (right-hand side of Figure \ref{fig:order}), the
error does not appreciably decrease at all (up to $2048$ degrees
of freedom per direction) for any of the numerical methods considered.
To run a simulation with $2048$ degrees of freedom in four dimensions
would require at least $280$ TB of memory, which would stress even
the world's most powerful supercomputers (the number one system on
the present TOP500 supercomputing list has $1310$ TB of memory; thus
such a simulation would be possible in principle, but certainly not
very practical). Thus, while the argument provided is valid for very
small times, it is all but pointless for realistic simulations. The
question that we should ask instead is how well does a numerical method
resolve the physics of the problem under consideration (for reasonable
problem sizes). 

The example considered here is two-dimensional and thus should not
be considered the gold standard for higher-dimensional simulations.
In fact, the reason why we consider this in two dimensions is that
it is prohibitively expensive to run simulations with $2048^{4}$
or even more degrees of freedom. This paper will show that for the
four-dimensional case the performance of a numerical method in the
asymptotic regime is markedly different from the performance of interest
in realistic simulations. More specifically, we will show that even
though the semi-Lagrangian discontinuous Galerkin method is inferior
to spline interpolation for small times, depending on the problem,
similar or even superior performance can be achieved for more realistic
simulations.

Of course, the importance of designing numerical methods that faithfully
resolve the physics (as opposed to primarily decreasing the error
in some norm) has been considered extensively in the context of long
time integration (see, for example, \cite{hairer2006geometric}).
The difference here is that we certainly would not consider $t=50$
for nonlinear Landau damping a long-time result. In fact, to determine
the long-time behavior of nonlinear Landau damping, numerical simulations,
mostly in two dimensions, for at least up to $t=1300$ have been considered
(see, for example, \cite{brunetti2000} and \cite{manfredi1997}).
In the literature much effort has been dedicated to numerical methods
which conserve certain invariants (such as mass or energy) exactly
and to methods which have good long time properties with respect to
these quantities. This is also an important aspect for the Vlasov
equation and many works exist that consider conservation of invariants
(see, for example, \cite{crouseilles2010,filbet2003} for spline interpolation
and \cite{einkemmer2017study,crouseilles2011} for the discontinuous
Galerkin approach). In addition, since the Vlasov equation is often
used to simulate plasma instabilities, the time evolution of certain
averaged quantities (most notably the total electric energy) are of
significant interest. Both of these aspects will be considered in
sections \ref{sec:nl} and \ref{sec:tsi}.

\section{Numerical methods\label{sec:numerical-method}}

To advance the numerical solution in time we use the splitting approach
introduced in \cite{cheng1976}. That is, we first consider the free-streaming
sub-flow
\begin{equation}
\partial_{t}f(t,x,v)+v\cdot\nabla_{x}f(t,x,v)=0,\qquad f(0,x,v)=g(x,v).\label{eq:free-streaming}
\end{equation}
The corresponding solution at time $\tau$ is then denoted by
\[
\mathrm{e}^{\tau A}g(x,v)=f(\tau,x,v).
\]
For the second part of the splitting algorithm we consider the acceleration
sub-flow
\begin{equation}
\partial_{t}f(t,x,v)+E(x)\cdot\nabla_{v}f(t,x,v)=0,\qquad f(0,x,v)=g(x,v),\label{eq:accelleration}
\end{equation}
where the electric field $E$ is determined from the (charge) density
$\rho(x)=\int g(x,v)\,\mathrm{d}v$ and is thus taken constant during
that step. This property actually follows from equation (\ref{eq:accelleration})
as the translation is only in the $v$-direction and thus leaves the
density invariant. We denote the corresponding solution by 
\[
\mathrm{e}^{\tau B}g(x,v)=f(\tau,x,v).
\]
Note that, contrary to what the notation suggests, $B$ is a nonlinear
operator as the electric field depends on $f$. Using this notation
we can easily formulate a time step of the splitting algorithm. In
all our numerical results we will use the second-order Strang splitting
given by
\[
f^{m+1}=\mathrm{e}^{\frac{\tau}{2}A}\mathrm{e}^{\tau B}\mathrm{e}^{\frac{\tau}{2}A}f^{m},
\]
where $f^{m}(x,v)$ is an approximation of $f(t_{m},x,v)$ and $\tau=t_{m+1}-t_{m}$
is the time step size. Let us note, however, that higher order splitting
methods have been constructed as well (see, for example, \cite{crouseilles2011high}). 

The main computational advantage of the splitting scheme is that the
resulting sub-flows, given by equations (\ref{eq:free-streaming})
and (\ref{eq:accelleration}), are in the form of an advection equation,
where the advection speed is independent of the variables being advected.
Thus, we can immediately solve equation (\ref{eq:free-streaming})
to obtain 
\[
\mathrm{e}^{\tau A}g(x,v)=g(x-\tau v,v)
\]
and solve equation (\ref{eq:accelleration}) to obtain 
\[
\mathrm{e}^{\tau B}g(x,v)=g(x,v-\tau E(x)).
\]
Thus, by using the splitting approach outlined above we have reduced
the task of computing a numerical approximation to the four-dimensional
Vlasov-{}-Poisson equation to computing two-dimensional translations
in phase space. However, this can be further simplified by splitting
equation (\ref{eq:free-streaming}) into an advection in the $x_{1}$
direction
\[
\partial_{t}f(t,x,v)+v_{1}\partial_{x_{1}}f(t,x,v)=0,\qquad f(0,x,v)=g(x,v)
\]
and an advection in the $x_{2}$ direction
\[
\partial_{t}f(t,x,v)+v_{2}\partial_{x_{2}}f(t,x,v)=0,\qquad f(0,x,v)=g(x,v).
\]
Let us denote the corresponding flows at time $\tau$ by $\mathrm{e}^{\tau A_{1}}g(x,v)$
and $\mathrm{e}^{\tau A_{2}}g$, respectively. Then since $[v_{1}\partial_{x_{1}},v_{2}\partial_{x_{2}}]=0$
(i.e.~the two operators commute) we have
\[
e^{\tau A_{2}}\mathrm{e}^{\tau A_{1}}g=e^{\tau A}g
\]
and no further splitting error is incurred. A very similar procedure
can be carried out for the acceleration sub-flow. Thus, the four-dimensional
Vlasov\textendash Poisson equation has been reduced to a sequence
of one-dimensional advections. This has significant advantages from
a computational point of view and, in particular, implies that we
only have to formulate our semi-Lagrangian scheme in a single dimension.
We should remark that this approach can be extended to more complicated
models, such as the Vlasov\textendash Maxwell equations \cite{crouseilles2015hamiltonian}.

We will now describe the two competing approaches considered in this
paper. First, for the semi-Lagrangian method based on spline interpolation
an equidistant grid is used. However, for any grid point $(x_{1i},x_{2j},v_{1k},v_{2l})$
the translations $x_{1i}-\tau v_{1k}$, $v_{1k}-\tau E_{1}(x_{1i},x_{2j})$,
etc., in general, do not coincide with a grid point (here $E_{1}$
denotes the first component of the electric field). Thus, we have
to use an interpolation scheme in order to evaluate $g(x_{1i}-\tau v_{1k},x_{2j},v_{1k},v_{2l})$,
$g(x_{1i},x_{2j},v_{1k}-\tau E_{1}(x_{1i},x_{2j}),v_{2l})$. This
is done by cubic spline interpolation. Note that only splines in one
space dimension have to be constructed. This is essential for the
computational efficiency of the scheme. Further computational ramifications
are discussed in more detail in section \ref{sec:performance}.

Second, the semi-Lagrangian discontinuous Galerkin scheme divides
the domain into cells $C_{ijkl}=I_{i}^{(1)}\times I_{j}^{(2)}\times I_{k}^{(3)}\times I_{l}^{(4)}$,
where the $I_{i}^{(1)},I_{j}^{(2)},I_{k}^{(3)},I_{l}^{(4)}$ are one-dimensional
intervals of length $h$. We further assume that a function $g$ is
given such that $g\vert_{C_{i}}$, i.e. the restriction of $g$ to
the $i$th cell, is a polynomial of degree $\ell$. Then the function
$g$ lies in the approximation space (note that we do not enforce
a continuity constraint across cell interfaces). However, in general,
this is not true for the translated function (for the other three
translations we proceed in exactly the same way)
\[
(T_{\tau}g)(x_{1},x_{2},v_{1},v_{2})=g(x_{1}-\tau v_{1},x_{2},v_{1},v_{2}).
\]
Thus, we perform an approximation by applying a projection operator
$P$ and obtain $PT_{\tau}g$. In order to turn this into a sensible
numerical scheme, it must hold that $Pf\vert_{C_{ijkl}}$ is a polynomial
of degree $\ell$ for any function $f$. Then $PT_{\tau}g$ constitutes
the sought-after approximation of $g(x_{1}-\tau v_{1},x_{2},v_{1},v_{2})$.
In the present case the operator $P$ is the $L^{2}$ projection on
the (finite dimensional) subspace of cell-wise polynomials of degree
$\ell$. Note that this projection is the only approximation made.
There is some freedom in selecting appropriate degrees of freedom
for this method. Two obvious choices are to use the coefficients in
a Legendre expansion or the value of function evaluations at the Gauss\textendash Legendre
points. In our implementation we use the latter. Since the value at
the beginning of each time step is (per assumption) a cell-wise polynomial
of degree $\ell$, we can use Gauss\textendash Legendre quadrature
in order to analytically compute the coefficients required to implement
the projection. The performance implications will be discussed in
section \ref{sec:performance}.

\section{Computational consideration \& performance results\label{sec:performance}}

The goal of this section is to discuss the performance of the cubic
spline and discontinuous Galerkin based semi-Lagrangian methods on
modern computer architectures. We will do this both from a theoretical
and a practical (i.e.~measuring specific implementations of these
algorithms) point of view.

As has been discussed in the previous section, only a set of one-dimensional
advection equations have to be solved for both algorithms. Thus, to
simplify the notation, we only consider the one-dimensional advection
\begin{equation}
u^{m+1}(x)=u^{m}(x-v\tau)\label{eq:1d-advection}
\end{equation}
in this section. This is, in fact, exactly how the numerical algorithm
is implemented on a computer. A function is provided that solves the
one-dimensional advection given by equation (\ref{eq:1d-advection}).
This function is then executed multiple times with different parameters.

For spline interpolation, the function values $u_{i}^{m}=u^{m}(x_{i})$
are given at the grid points $x_{i}$. The index $i$ runs from $0$
to $n-1$, where $n$ is the number of degrees of freedom. From this
a spline that covers the entire domain has to be constructed. Since
it would be wasteful to represent a spline as a collection of piecewise
polynomials (as most of the coefficients of these polynomials are
fixed by the continuity constraints), a B-spline basis is used instead.
That is, the spline interpolation $\tilde{u}^{m}$ of $u_{i}^{m}$
is represented as
\[
\tilde{u}^{m}(x)=\sum_{k}\omega_{k}S(x-x_{k}),
\]
where $S$ is given by
\[
6S(x)=\begin{cases}
4-6(x/h)^{2}+3\vert x/h\vert^{3} & 0\leq\vert x\vert\leq h\\
(2-\vert x/h\vert)^{3} & h\leq\vert x\vert\leq2h\\
0 & \text{otherwise}
\end{cases}
\]
and $h$ is the grid spacing. The coefficients $\omega_{k}$ are uniquely
determined by the function values and can be obtained by solving the
following linear system of equations (assuming periodic boundary conditions)

\[
\frac{1}{6}\left[\begin{array}{cccccc}
4 & 1 & 0 & \cdots & 0 & 1\\
1 & 4 & 1 & 0 & \cdots & 0\\
0 & 1 & 4 & 1 & \ddots & \vdots\\
0 & 0 & \ddots & \ddots & \ddots & 0\\
\vdots & \vdots & \ddots & 1 & 4 & 1\\
1 & 0 & \cdots & 0 & 1 & 4
\end{array}\right]\left[\begin{array}{c}
\omega_{0}\\
\omega_{1}\\
\vdots\\
\vdots\\
\vdots\\
\omega_{n-1}
\end{array}\right]=\left[\begin{array}{c}
u_{0}\\
u_{1}\\
\vdots\\
\vdots\\
\vdots\\
u_{n-1}
\end{array}\right].
\]
Note that the matrix is tridiagonal and symmetric and thus efficient
methods exist to solve this system. Once the spline is constructed
it is used to solve equation (\ref{eq:1d-advection}) as follows
\[
u_{i}^{m+1}=\tilde{u}^{m}(x_{i}-v\tau).
\]
Thus, the numerical algorithm proceeds in two steps. First, the spline
is constructed by solving a linear system and this intermediate representation
is then used to evaluate the numerical solution at the grid points.
Both of these steps require $\mathcal{O}(n)$ memory and arithmetic
operations, where $n$ is the number of degrees of freedom. In the
construction step this is due to the fact that linear complexity solvers
are available for tridiagonal systems (such as the Thomas algorithm)
and for the evaluation step this is true since the support of $S$
is compact and independent of $n$. 

On modern computer hardware both of these steps are memory bound.
That is, performance is dictated by how fast data can be read and
written to memory (as opposed to how many arithmetic operations the
computer can perform). In the present case we have to at least read
all $u_{i}^{m}$ and write all $u_{i}^{m+1}$. In addition, we have
to write and read all $\omega_{k}$ at least once. The latter is more
subtle than one might expect. Since in a sequential implementation
only one one-dimensional array of $\omega_{k}$s needs to exist at
any given time, this array could conceivable be kept in cache (except
for extremely large problem sizes) and then reads and writes to it
would not count as a memory operation (in the sense that those do
not negatively impact performance).

\begin{figure}
\centering{}\includegraphics[width=6cm]{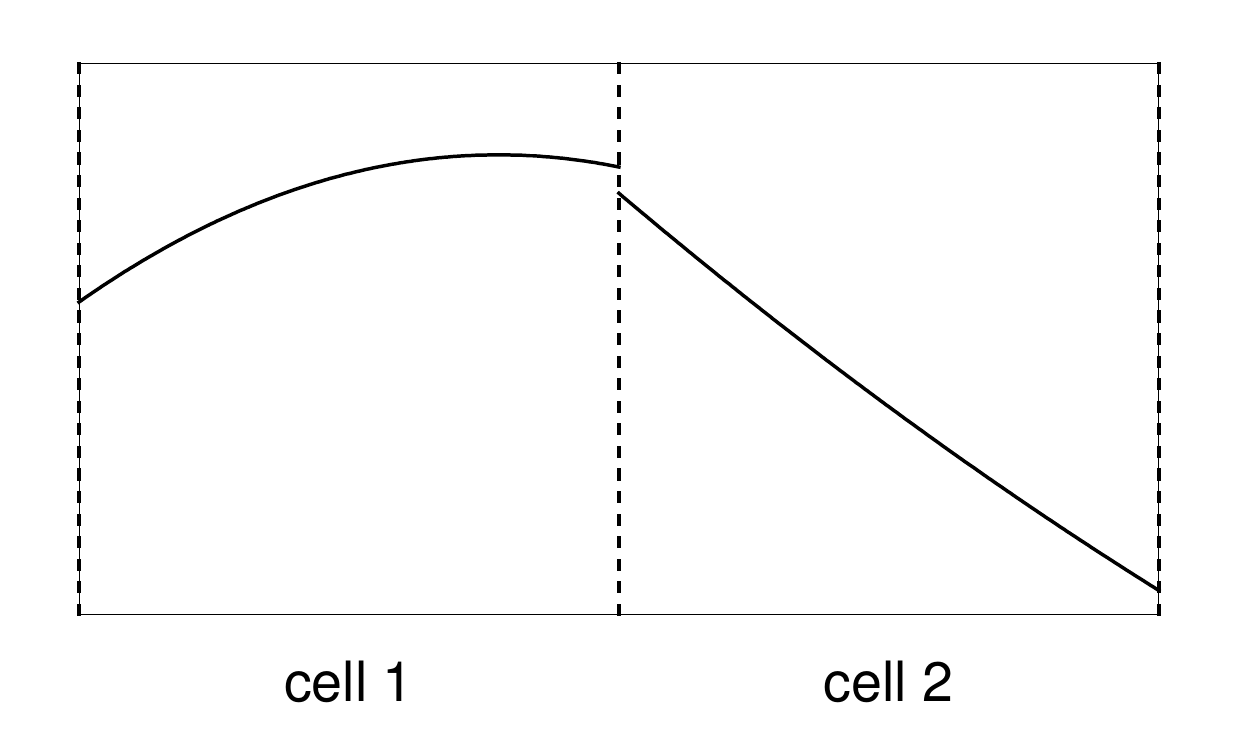}\includegraphics[width=6cm]{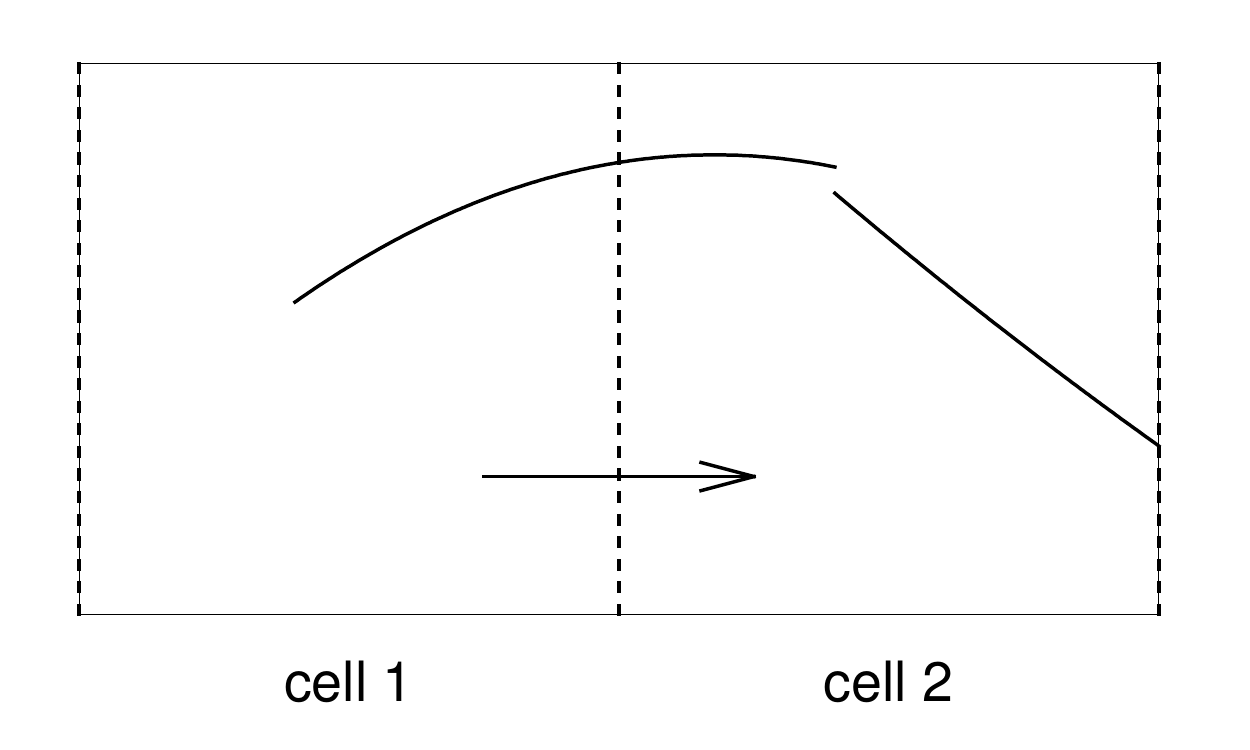}\caption{Illustration of the semi-Lagrangian discontinuous Galerkin scheme.
The piecewise polynomial function of degree $\ell$ (left picture)
is translated by $v\tau$ (right picture). The translated function
is then projected back to the approximation subspace (which requires
data from at most two adjacent cells).\label{fig:illustration-sLdG}}
\end{figure}

For the semi-Lagrangian discontinuous Galerkin scheme the degrees
of freedom are given by function values $u_{ij}^{m}\approx u^{m}(x_{i-1/2}+\xi_{j})$,
where $i$ is the cell-index (which runs from $0$ to $n_{C}$), $x_{i-1/2}$
the left cell interface, and $\xi_{j}$ is the $j$th Gauss\textendash Legendre
node scaled to the interval $[0,h]$. The indices $j$ run from $0$
to $\ell$ (the degree of the polynomial approximant). The resulting
scheme has order $\ell+1$ and we will henceforth denote it as dG($\ell+1$).
It is clear that to perform the $L^{2}$ projection on the $i$th
cell only data from at most two adjacent cells is required. This is
illustrated in Figure \ref{fig:illustration-sLdG}. We denote the
index of these two cells by $i^{\star}$ and $i^{\star}+1$, respectively.
Since the projection is a linear operator we can write the resulting
numerical scheme as follows

\begin{equation}
u_{ij}^{m+1}=\sum_{l}A_{jl}u_{i^{\star}l}^{m}+\sum_{j}B_{jl}u_{i^{\star}+1;l}^{m},\label{eq:AB-sldg}
\end{equation}
where $A\in\mathbb{R}^{(\ell+1)\times(\ell+1)}$ and $B\in\mathbb{R}^{(\ell+1)\times(\ell+1)}$
are matrices that only depend on $v\tau$. The matrices can be precomputed
at the beginning of each time step and thus their cost will not be
a major concern except for extremely coarse space discretizations.
The derivation of their exact form is somewhat lengthy and we thus
refer the interested reader to \cite{einkemmer2015}.

For the implementation equation (\ref{eq:AB-sldg}) is employed. The
size of the matrices is negligible in comparison to the cache size
and we therefore do not have to count reads and writes to them as
memory accesses. The present algorithm has to read all entries in
$u_{ij}^{m}$ and write to all entries in $u_{ij}^{m+1}$ ; it thus
requires $\mathcal{O}(n)$ memory operations, where $n=n_{C}(\ell+1)$
is the number of degrees of freedom. Technically, the algorithm requires
$\mathcal{O}(\ell n)$ arithmetic operations and thus does not scale
linearly in the degrees of freedom. However, since, for reasonable
values of $\ell$, this is still an extremely memory bound problem
this is of no consequence (i.e.~performance is dictated by the memory
accesses and not by amount of arithmetic operations that need to be
performed).

In an idealized setting both algorithms could achieve very similar
levels of performance (both require us to read $u^{m}$ once and write
to $u^{m+1}$ once). Of course, achieving this in an actual implementation
is not a trivial task. However, it is our believe that the local nature
of the discontinuous Galerkin approach simplifies this task, at least,
to some extent. In the following we will compare the performance of
two Vlasov solvers. SeLaLib \cite{selalib} is a library that facilitate
kinetic simulations based on the semi-Lagrangian approach. It has
been used extensively (see, for example, \cite{crouseilles2014new,guisset2016lagrangian,barsamian2017verification}).
The library is written in Fortran and provides an efficient implementation
of the cubic spline based semi-Lagrangian discontinuous Galerkin method.
SeLaLib is compared with our own implementation (called sldg) of the
semi-Lagrangian discontinuous Galerkin approach. The code is written
in C++ and due to the use of templates facilitates implementations
that work in arbitrary dimensions, while still maintaining good performance.
This architecture is described in \cite{einkemmer2015}. In addition
to solving the Vlasov equation, it has been used in non-traditional
applications of semi-Lagrangian methods (such as for the Kadomtsev\textendash Petviashvili
equation that requires the solution of a nonlinear advection \cite{einkemmer2017semi}).

Another major difference in the implementation of these two codes
is how data is presented to the routines that perform the one-dimensional
advections. For SeLaLib the corresponding arrays are always contiguous
in memory. That is, for the advection in the $x_{1}$-direction the
elements in the $x_{1}$ direction are laid out sequentially in memory,
for the advection in the $x_{2}$-direction the elements in the $x_{2}$
direction are laid out sequentially in memory, etc. This requires
that the four dimensional array is converted (i.e.~transposed) to
the appropriate form for each advection. This approach has the advantage
that the same code can be used for all advection steps and cache locality
is ensured almost automatically. It also fits very well with spline
interpolation, as the global algorithm for constructing the spline
is pointed at one contiguous array of memory. The disadvantage of
this approach is that the transpositions do not have a very favorable
memory access pattern and thus incur a significant computational cost.

On the other hand, the semi-Lagrangian discontinuous Galerkin implementation
uses the ability of abstractions in C++ to hide the actual memory
layout from the function that performs the advection. Thus, it is
still possible to write the one-dimensional advection code as if the
arrays are contiguous, even though the actual memory layout is more
complicated. As a consequence, no transpositions have to be carried
out in this implementation. The main concern here is that cache locality
has to be taken into account. However, the local nature of the discontinuous
Galerkin approach makes this less problematic. To further alleviate
this issue we use cache blocking for the dimension with the largest
stride.

We have measured both the performance and the amount of memory used
for these two codes. The results are shown in Table \ref{tab:runtime-mem-comparison}.
For both methods one dual socket node of the described computer systems
is used. Both codes are run using $16$ parallel threads. We observe
that the discontinuous Galerkin implementation is both faster (by
approximately a factor of $2$ on the newer Haswell CPUs and approximately
a factor of $1.7$ on the Ivy Bridge CPUs) and uses significantly
less memory (by approximately a factor of $3$), if the same number
of degrees of freedom are used.

\begin{table}
\centering
\textbf{Workstation}
\vspace{0.1cm}

\begin{tabular}{rrrrrr}
& \multicolumn{2}{c}{Time per step [s]} & & \multicolumn{2}{c}{Memory used [GB]}\\
\cline{2-3} \cline{5-6}
& dof=64 & dof=128 & & dof=64 & dof=128 \\
\hline
SeLaLib     & 0.194  & 3.25  & & 0.85 & 13.0 \\
sldg (dG4)  & 0.105  & 1.62  & & 0.28 & 4.2 \\
\hline
\end{tabular}
\vspace{0.3cm}

\textbf{VSC3}
\vspace{0.1cm}

\begin{tabular}{rrr}
& \multicolumn{2}{c}{Time per step [s]}\\
\cline{2-3}
& dof=64 & dof=128 \\
\hline
SeLaLib    & 0.219 & 4.02 \\
sldg (dG4) & 0.151 & 2.31 \\
\hline
\end{tabular}

\caption{Comparison of performance and memory usage between SeLaLib (cubic
spline interpolation) and sldg (fourth order semi-Lagrangian discontinuous
Galerkin method). For both codes the nonlinear Landau damping initial
value of section \ref{sec:nl} is used and $64$ and $128$ degrees
of freedom (dof) per direction are considered. The Workstation consists
of two Intel Xeon E5-2630 v3 CPUs and 32 GB of DDR4 memory. One node
of the VSC3 consists of two Intel Xeon E5-2650v2 and 64GB of DDR3
memory. The memory consumption is measured using the \texttt{/usr/bin/time}
utility. \label{tab:runtime-mem-comparison}}
\end{table}
Both codes support parallelization using MPI. This is essential even
to obtain good performance on a single node (i.e.~for the results
presented in Table \ref{tab:runtime-mem-comparison}). However, both
codes can also be run on distributed memory clusters (supercomputers).
In this case the issue of how to distribute the computational domain
to the different MPI processes becomes important. Although this choice
of data distribution is identical for both codes it is instructive
to discuss it. The main issue here is that most FFT libraries (in
particular, FFTW, which is used in both codes) only allow a so-called
pencil decomposition. Thus it is not straightforward to parallelize
the Fourier transforms (without incurring a penalty associated with
additional data transfer). Both codes avoid this problem by only parallelizing
in the velocity directions. Then, the Fourier transforms can be conducted
independently on each of the processes. For four/six dimensional simulations
there are still two/three dimensions available for parallelization.
The main difference in the two codes is that while sldg only performs
halo communication (due to the local nature of the semi-Lagrangian
discontinuous Galerkin scheme), SeLaLib requires two calls to \texttt{MPI\_Alltoall
}in each time step (due to the global data dependency inherent in
constructing the spline).

\begin{figure}
\begin{centering}
\includegraphics[width=6cm]{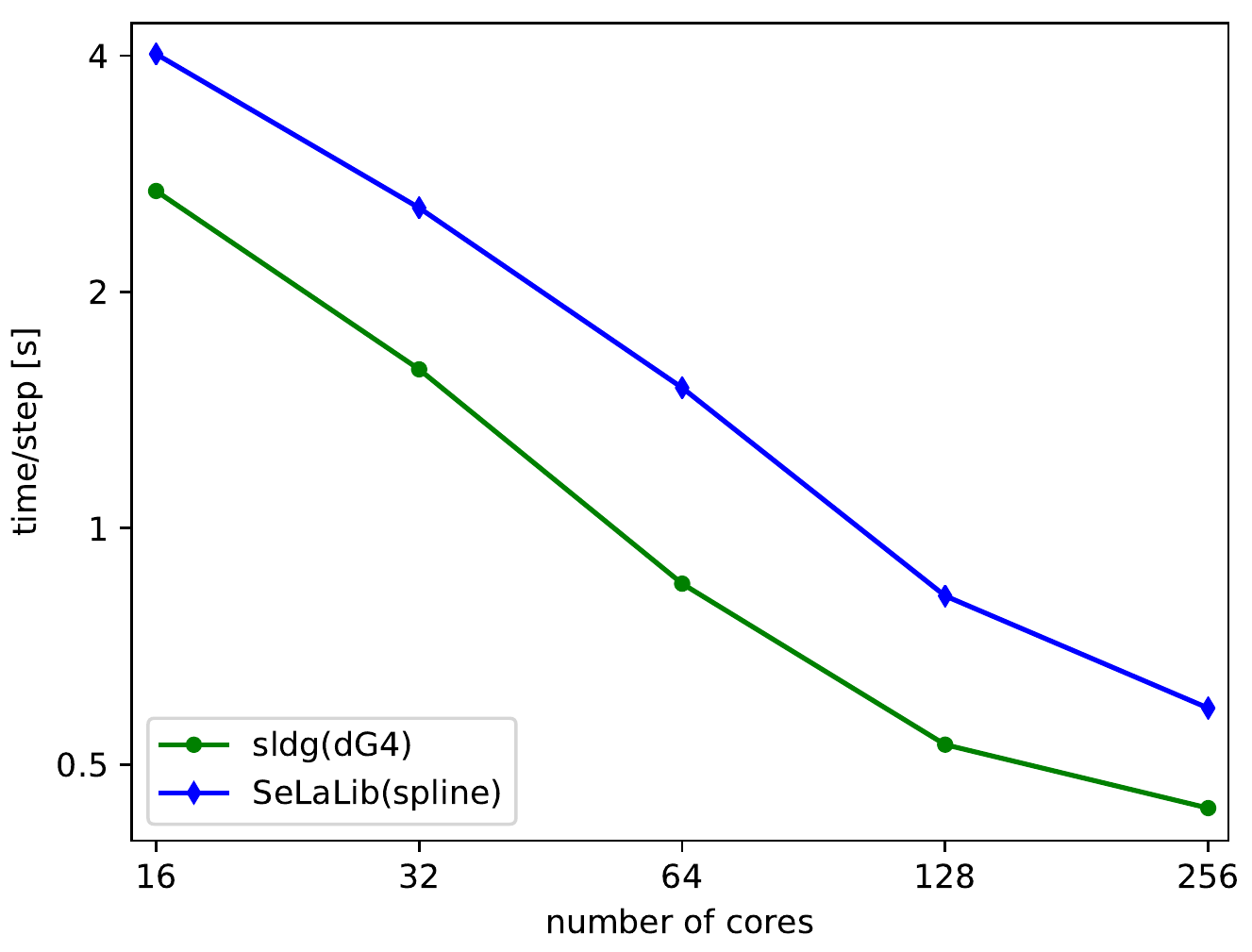}\includegraphics[width=6cm]{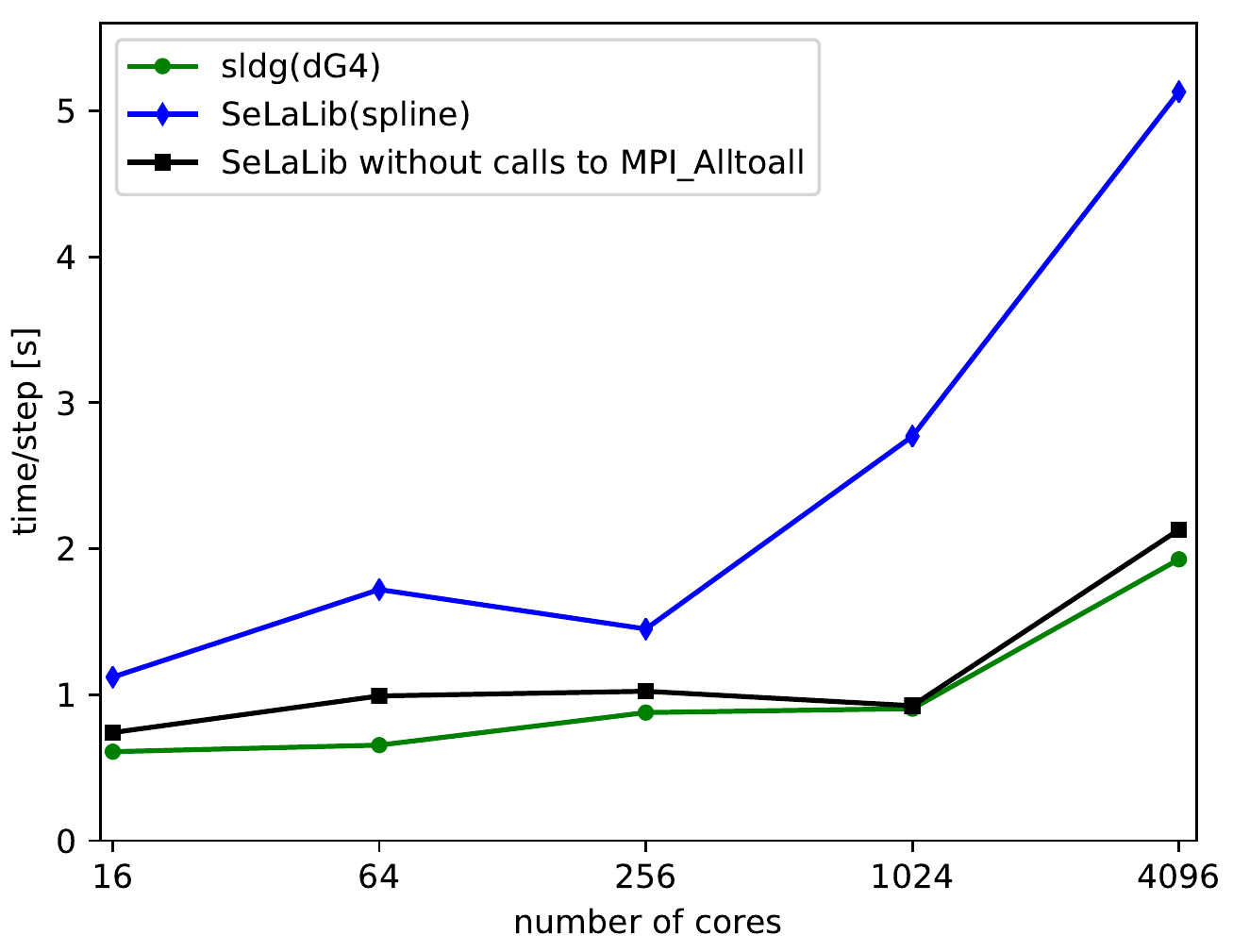}
\par\end{centering}
\caption{Strong and weak scaling for SeLaLib (cubic spline interpolation) and
sldg (fourth order semi-Lagrangian discontinuous Galerkin method).
For strong scaling the nonlinear Landau damping initial value of section
\ref{sec:nl} with $128$ degrees of freedom per direction is used.
For weak scaling $64^{2}\times(\sqrt{n}\cdot128)^{2}$ degrees of
freedom, where $n$ is the number of nodes employed, are used. The
simulations are conducted on the VSC3, each node of which consists
of two Intel Xeon E5-2650v2 and 64GB of DDR3 memory.\label{fig:scaling}}
\end{figure}

We now investigate the scaling of both codes on the Vienna Scientific
Cluster 3 (VSC3). For more information and details on the hardware
of that system we refer the reader to \url{http://vsc.ac.at/systems/vsc-3/}.
On the left hand side of Figure \ref{fig:scaling} we present strong
scaling results (i.e.~the problem size is kept constant as the number
of cores is increased) for the problem considered in section \ref{sec:nl}.
We observe that the discontinuous Galerkin scheme maintains its advantage
independently of the number of cores the program is run on (i.e.~both
codes show a similar behavior with respect to scaling). On the right
hand side of Figure \ref{fig:scaling} we present weak scaling results
(i.e.~the problem size increases proportionally to the number of
cores). In this case we observe that sldg significantly outperforms
SeLaLib (in particular, as we consider more than $256$ cores). It
is interesting to investigate the reason for this discrepancy in performance.
To that end we have subtracted the time required for the calls to
\texttt{MPI\_Alltoall,} which are used to map the velocity grid into
contiguous memory on each process, in the SeLaLib code from the total
run time (note that even on a single node this operation requires
a significant amount of time, as computing the transpose is a memory
intensive operation). The resulting line is very similar to the run
time of the sldg code. That is, almost all of the difference in performance
between the two codes is explained by the global communication required
for the spline interpolation. 

\textbf{}

\section{Nonlinear Landau damping\label{sec:nl}}

In this section we consider nonlinear Landau damping in four dimensions.
That is, we impose the initial value
\begin{equation}
f(0,x_{1},x_{2},v_{1},v_{2})=\frac{1}{2\pi}\left(1+\epsilon(\cos(kx_{1})+\cos(kx_{2}))\right)\mathrm{e}^{-(v_{1}^{2}+v_{2}^{2})/2}\label{eq:iv-nl}
\end{equation}
on the domain $[0,4\pi]^{2}\times[-6,6]^{2}$. We set $\epsilon=\frac{1}{2}$
and $k=0.5$. This is precisely the setup used in \cite{filbet2003}.
Strang splitting with step size $\tau=0.1$ is used for the time integration.
We are interested in resolving the dynamics of nonlinear Landau damping
up to $t=50$. More precisely, the characteristic growth and decay
of the electric energy should be indistinguishable (in a plot) from
the exact solution until $t=50$. To accomplish this for spline interpolation
we need $128$ grid points per direction (see Figure \ref{fig:order}).
Note that despite the fact that, as the results in Figure \ref{fig:order}
show, the error in the density function is on the order of unity,
the physically important quantity is resolved very well. Therefore,
there is simply no advantage of running the integrator in the asymptotic
regime (which in any case would not be feasible in four dimensions).

In the introduction second, fourth, and sixth order discontinuous
Galerkin schemes were considered. It was already clear at that point
that the second order scheme is not competitive. However, the order
plots (for $t=10$) presented there, seem to indicate that the sixth
order discontinuous Galerkin approach performs best. However, this
is misleading as numerical simulations have shown that for $t=50$
(where the plots in the introduction give no sensible indication)
the fourth order method should actually be favored. Therefore, in
this section we will exclusively use dG4. Before proceeding, let us
note, however, that for long time integration this is not necessarily
the case. For example \cite{einkemmer2017study} found by performing
numerical simulations in two dimensions that for long times (both
for nonlinear Landau damping and for a bump on tail instabilities)
the sixth order method can be significantly more efficient.

If we run the same simulation for the discontinuous Galerkin approach
we see that the solution is somewhat worse. In particular, between
$t=40$ and $t=50$ the growth rate of the electric energy is slightly
increased compared to the reference solution. In order to remedy this
we can increase the number of degrees of freedom per direction to
$144$ (corresponding to $\alpha=1.125$). This means that in four
dimensions 60\% more degrees of freedom are required by the discontinuous
Galerkin approach to achieve the same accuracy. However, as we have
seen in section \ref{sec:performance}, the implementation of the
discontinuous Galerkin scheme is also more efficient (by a factor
of $2$ on Haswell CPUs) and uses less memory (by a factor of $3$).
Thus, the discontinuous Galerkin scheme retains an advantage in performance
(approximately 25\%). In addition, the amount of memory required is
reduced by a factor of $2$.

\begin{figure}[H]
\begin{centering}
\includegraphics[width=6cm]{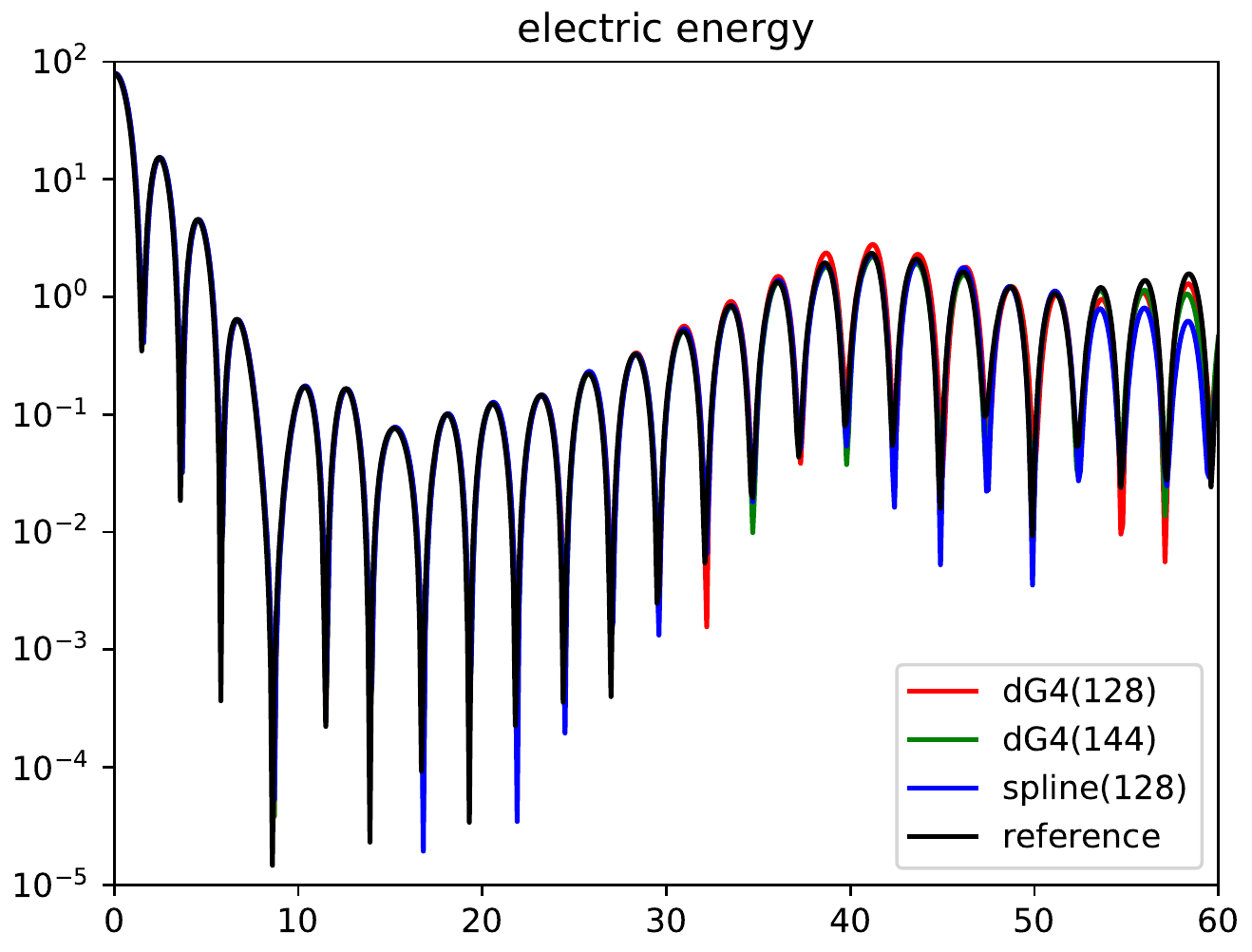}\includegraphics[width=6cm]{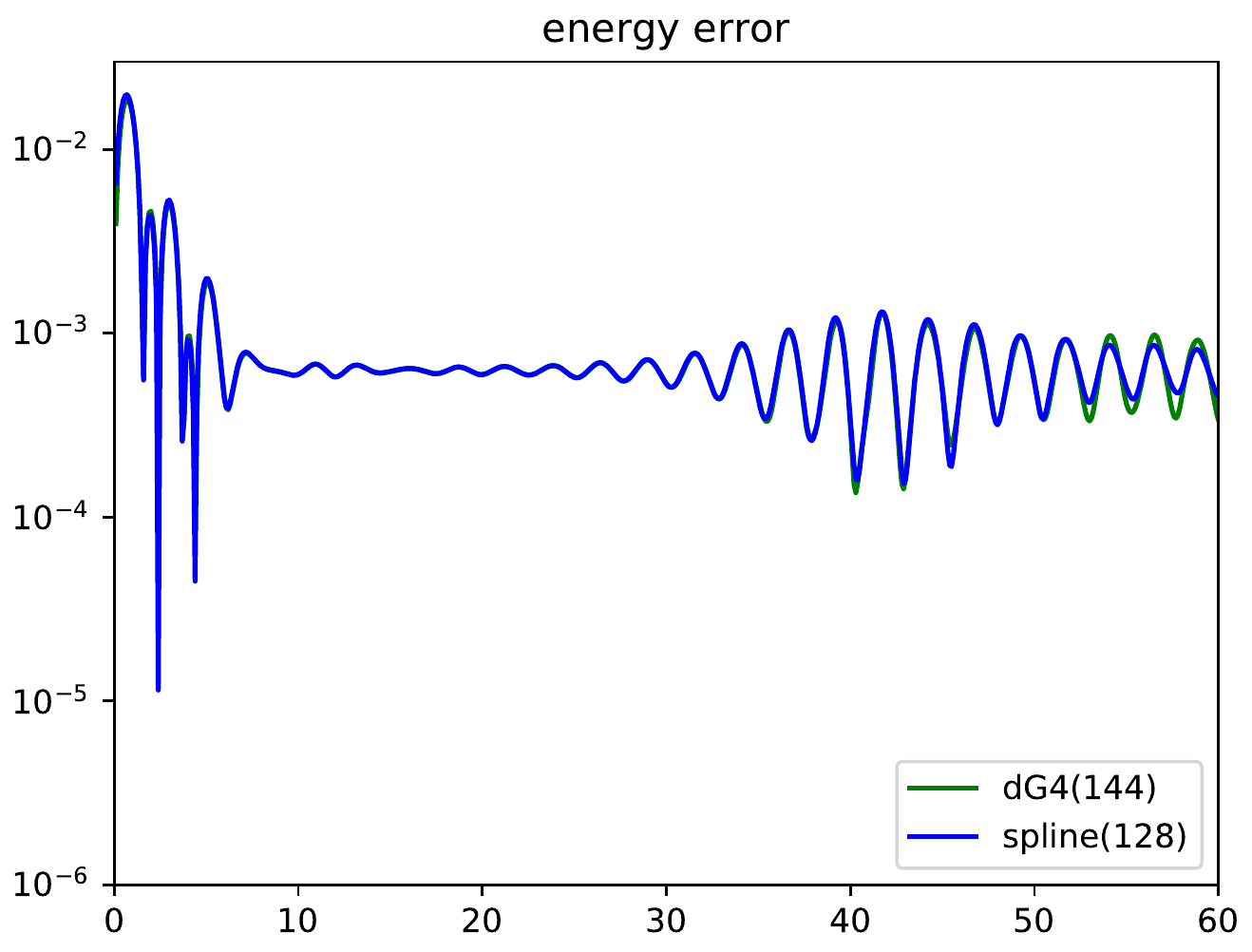}
\par\end{centering}
\begin{centering}
\includegraphics[width=6cm]{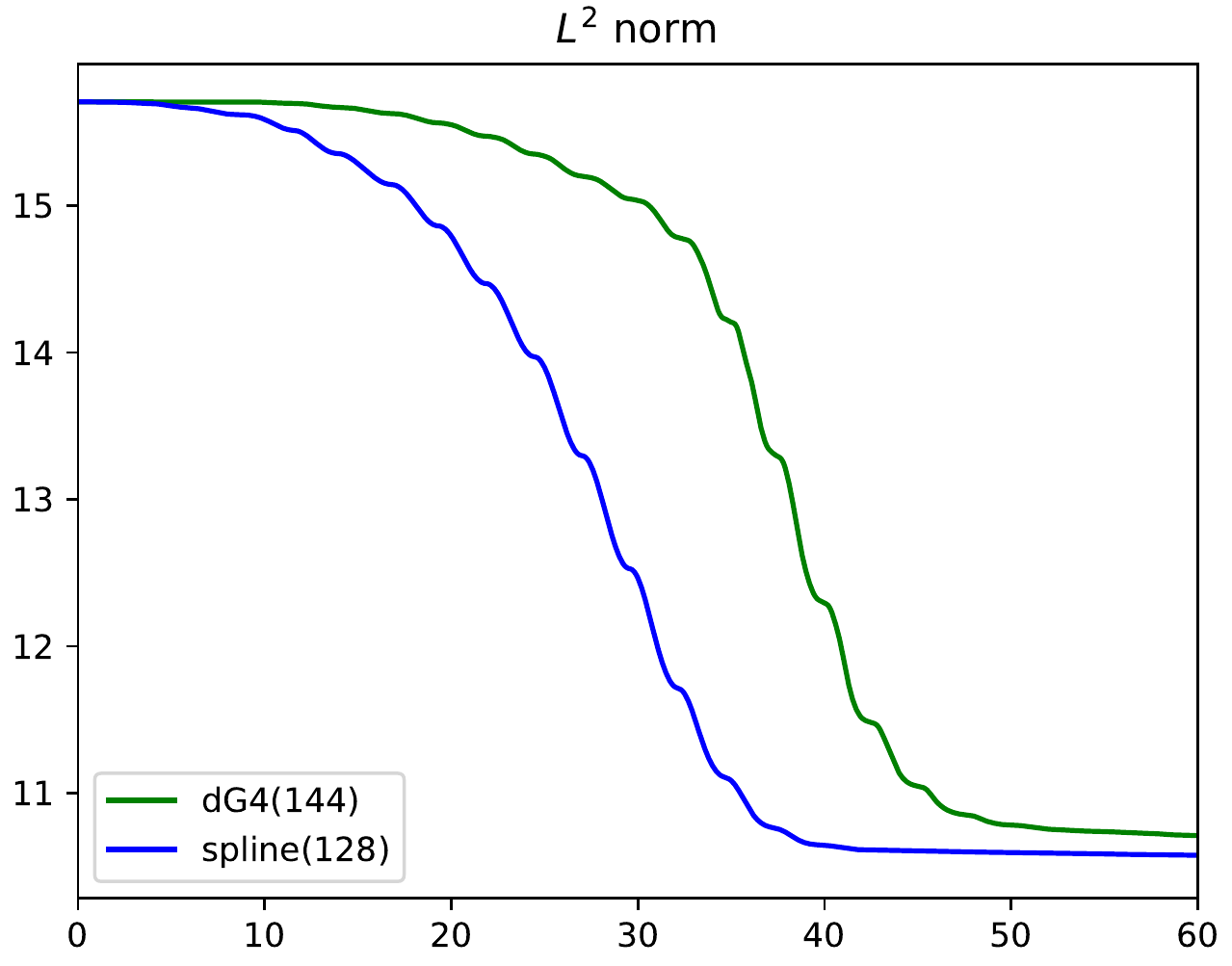}\includegraphics[width=6cm]{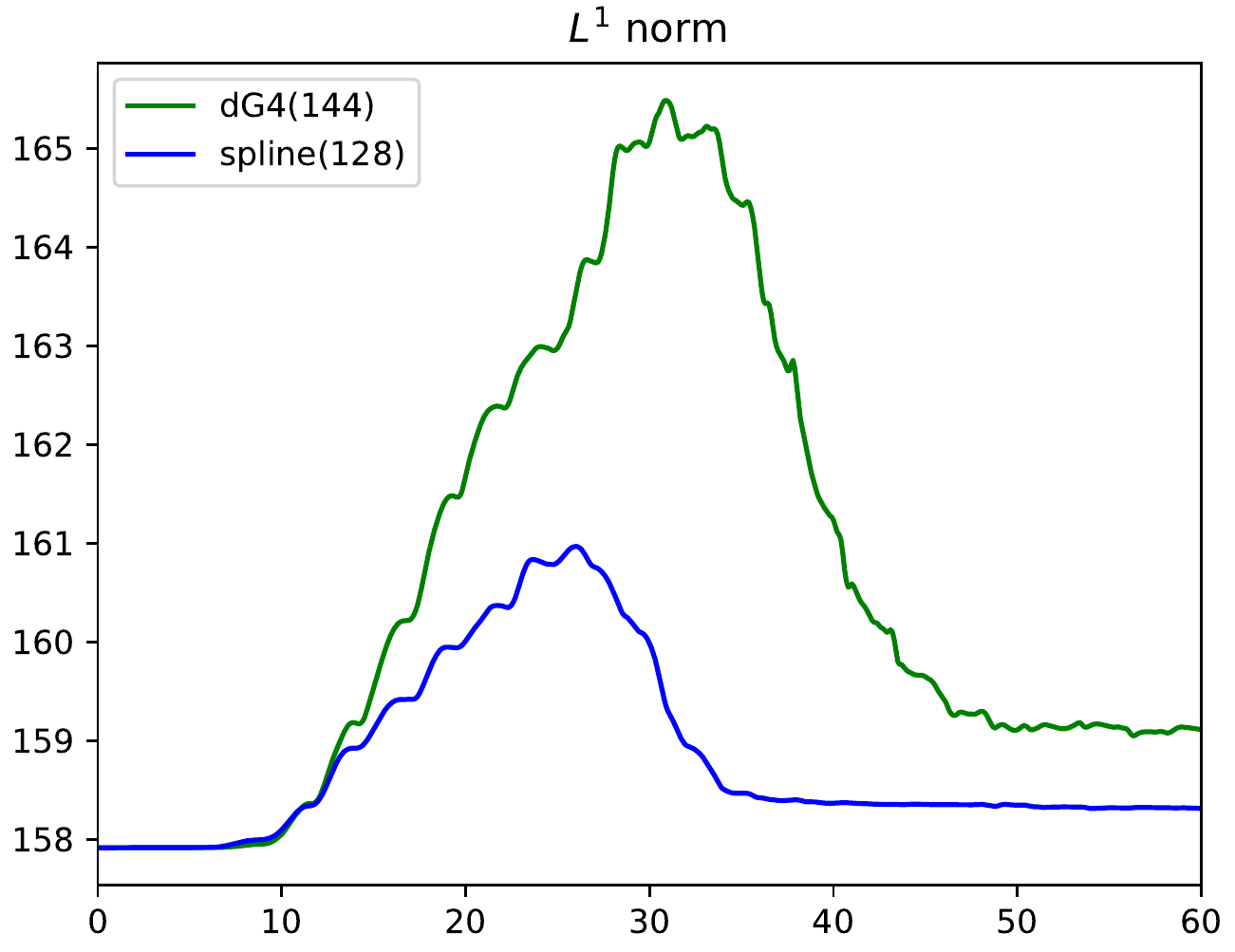}
\par\end{centering}
\caption{The electric energy (top left), the energy error (top right), the
$L^{2}$ norm (bottom left), and the $L^{1}$ norm (bottom right)
are shown as a function of time for the nonlinear Landau damping problem
given in equation (\ref{eq:iv-nl}). The reference solution is computed
using spline interpolation and $512$ grid points (per direction).
\label{fig:nl}}
\end{figure}

In addition, Figure $\ref{fig:nl}$ shows the important physical invariants
of energy, $L^{2}$ norm, and $L^{1}$ norm. Although the main concern
for the nonlinear Landau damping is the time evolution of the electric
energy, these invariants give additional diagnostics of how well physical
properties of the exact solution are respected by a numerical method.

The error in energy is almost identical for both methods and is on
the order of $10^{-3}$ (after a transient phase). This is due to
the fact that the energy error for the present problem is dominated
by the time integration error.

For collisionless plasma simulations the $L^{2}$ error is an important
diagnostic. A decrease in the $L^{2}$ norm implies that our numerical
scheme adds numerical diffusion to the simulation. Although most numerical
methods introduce a certain amount of diffusion, it is clearly preferable
to limit this as much as possible. For nonlinear Landau damping the
discontinuous Galerkin method outperforms spline interpolation in
this metric (in particular, between $t=10$ and $t=40$).

Finally, we look at the $L^{1}$ norm. For the exact solution mass
and $L^{1}$ norm are equal in value and the $L^{1}$ norm is thus
conserved. However, if a mass conservative numerical scheme commits
a numerical error that introduces negative values, the result is a
change in the $L^{1}$ norm. This is, in particular, a problem for
schemes that are susceptible to oscillations (such as spectral methods).
However, most numerical methods for the Vlasov equation experience
this behavior to some extent. Let us emphasize that (contrary to say
the Euler equations) this neither leads to numerical instabilities
nor necessarily diminishes averaged quantities (such as the electric
energy). However, since negative values make the physical interpretation
of the density function difficult, it is still desirable to limit
this behavior as much as possible. For the nonlinear Landau damping
the discontinuous Galerkin scheme shows worse performance with respect
to $L^{1}$ conservation compared to spline interpolation.

Overall we conclude that the discontinuous Galerkin implementation
has a moderate advantage with respect to performance and $L^{2}$
conservation as well as a significant advantage with respect to memory
consumption. On the other hand, spline interpolation has a moderate
advantage with respect to $L^{1}$ conservation. If we assume that
the present results generalize to the six-dimensional case (i.e.~$\alpha=1.125$
as above), then the performance of both methods would be evenly matched.
This is due to the fact that $1.125^{6}\approx2$ more degrees of
freedom would be required for the discontinuous Galerkin scheme, which
exactly cancels the factor of $2$ advantage the corresponding implementation
provides. 

\section{Two-stream instability\label{sec:tsi}}

In this section we consider the two-stream instability in four dimensions.
That is, we impose the initial value
\begin{align}
f(0,x_{1},x_{2},v_{1},v_{2}) & =(1+\epsilon\cos(kx_{1})\cos(kx_{2}))f^{\text{eq}}(v_{1},v_{2})\label{eq:iv-tsi}
\end{align}
with 
\[
f^{\text{eq}}=\frac{1}{8\pi}\left(\mathrm{e}^{-(v_{1}-v_{0})^{2}/2}+\mathrm{e}^{-(v_{1}+v_{0})^{2}/2}\right)\left(\mathrm{e}^{-(v_{2}-v_{0})^{2}/2}+\mathrm{e}^{-(v_{2}+v_{0})^{2}/2}\right),
\]
where $\epsilon=10^{-3}$, $k=0.2$, and $v_{0}=2.4$. The computational
domain is $[0,10\pi]^{2}\times[-6,6]^{2}$. Strang splitting with
step size $\tau=0.1$ is used for the time integration. 

This is an interesting problem as the electric energy is small (even
decreases by three orders of magnitude compared to the initial value)
during a quiescent phase (from $t=0$ to approximately $t=175$).
However, eventually a physical instability develops that results in
an exponential increase of the electric energy (approximately $t=175$
to $t=225$). This is followed by nonlinear saturation. In this nonlinear
phase the electric energy oscillates around a constant value with
an amplitude that decreases in time. 

It turns out that to obtain agreement between a numerical and the
exact solution for the quiescent phase is challenging. In fact, simulations
conducted for a similar two-dimensional problem indicate that for
spline interpolation this is only achieved for $1024$ grid points
per direction (which is extremely expensive). Nevertheless, simulations
with fewer degrees of freedom are still able to give pertinent information
on the plasma instability. In particular, at which time the instability
occurs, the growth rate of the instability, and the behavior of the
electric energy in the nonlinear regime (i.e. after saturation).

\begin{figure}
\begin{centering}
\includegraphics[width=6cm]{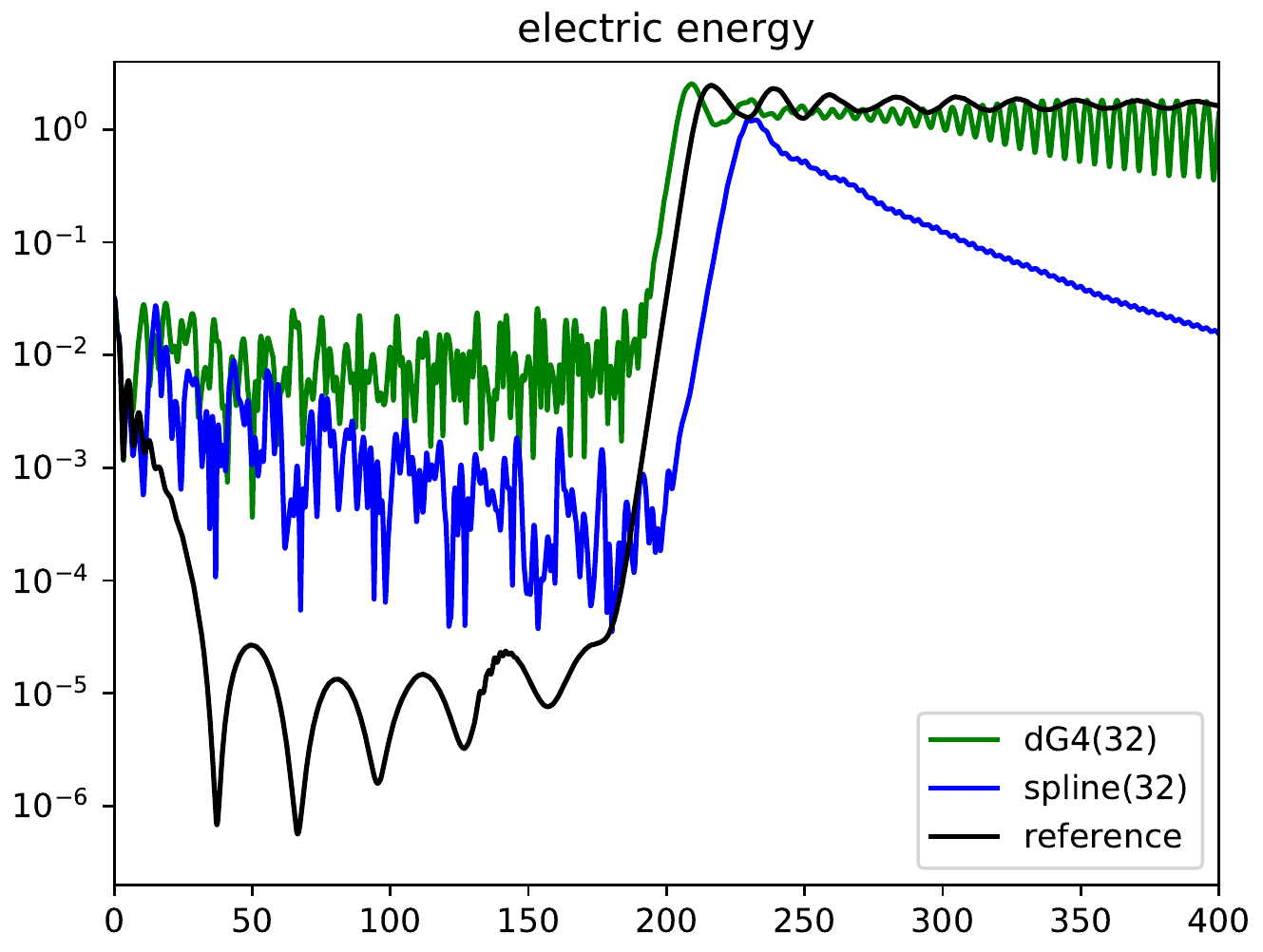}\includegraphics[width=6cm]{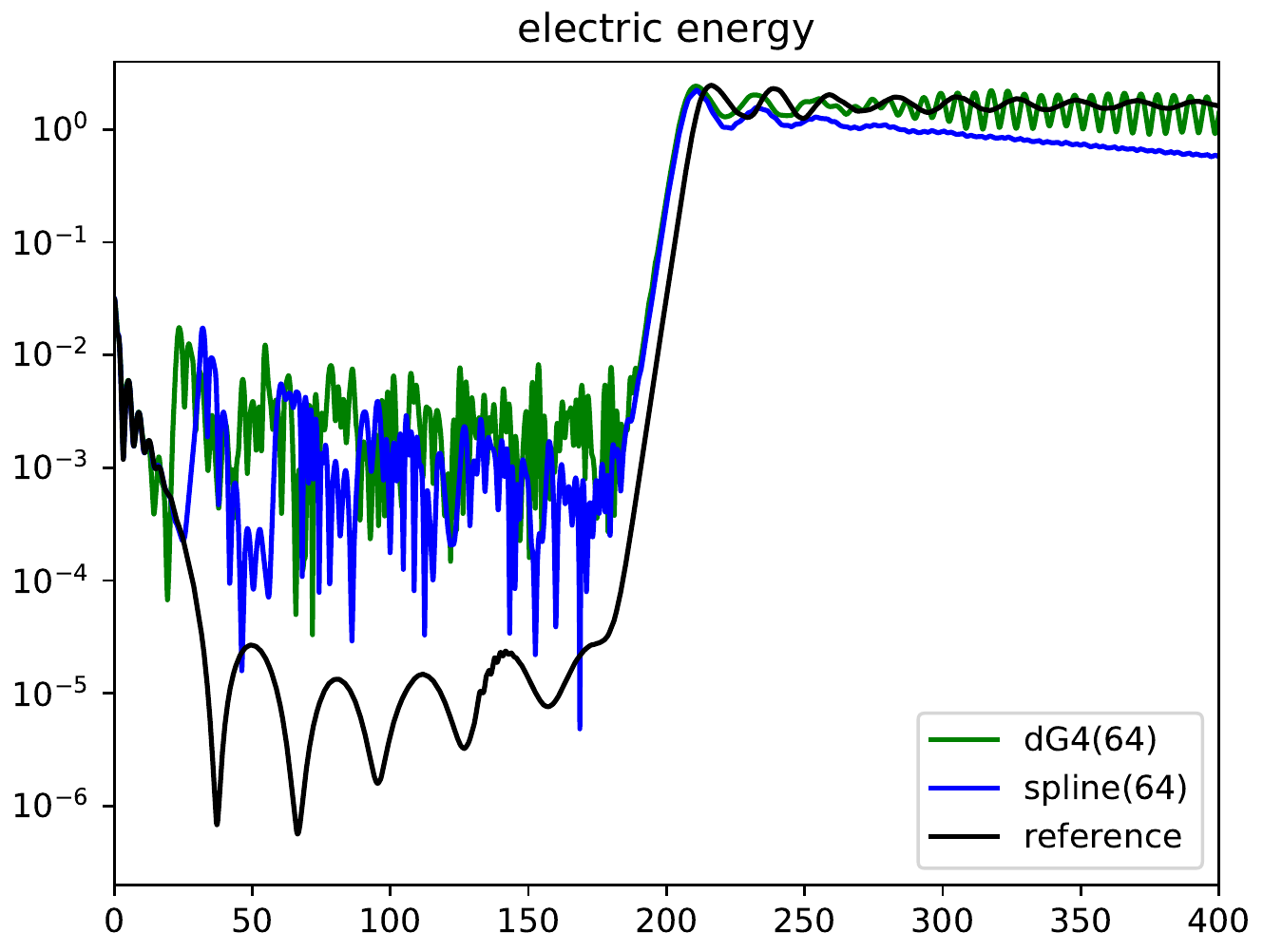}
\par\end{centering}
\begin{centering}
\includegraphics[width=6cm]{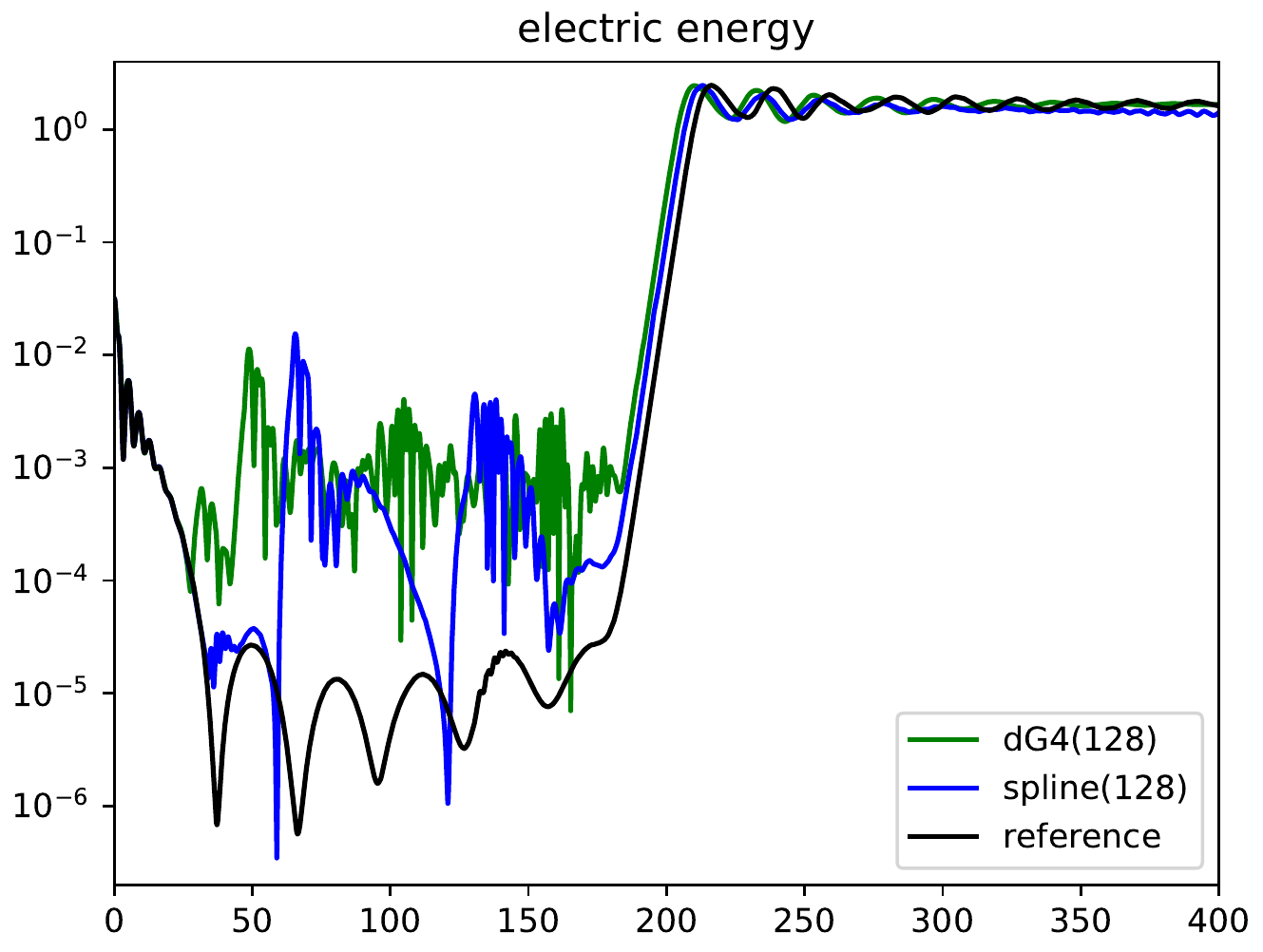}\includegraphics[width=6cm]{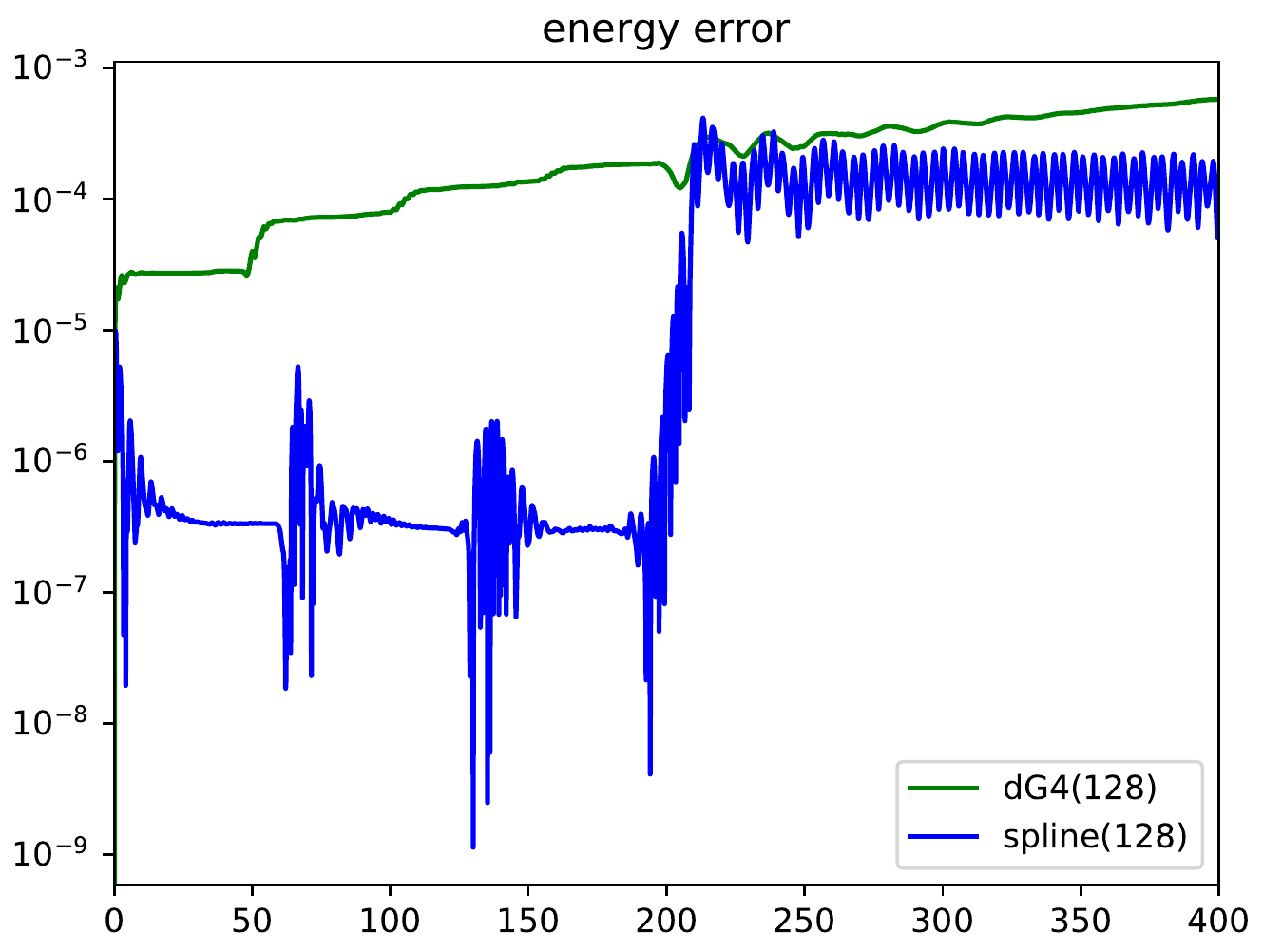}
\par\end{centering}
\begin{centering}
\includegraphics[width=6cm]{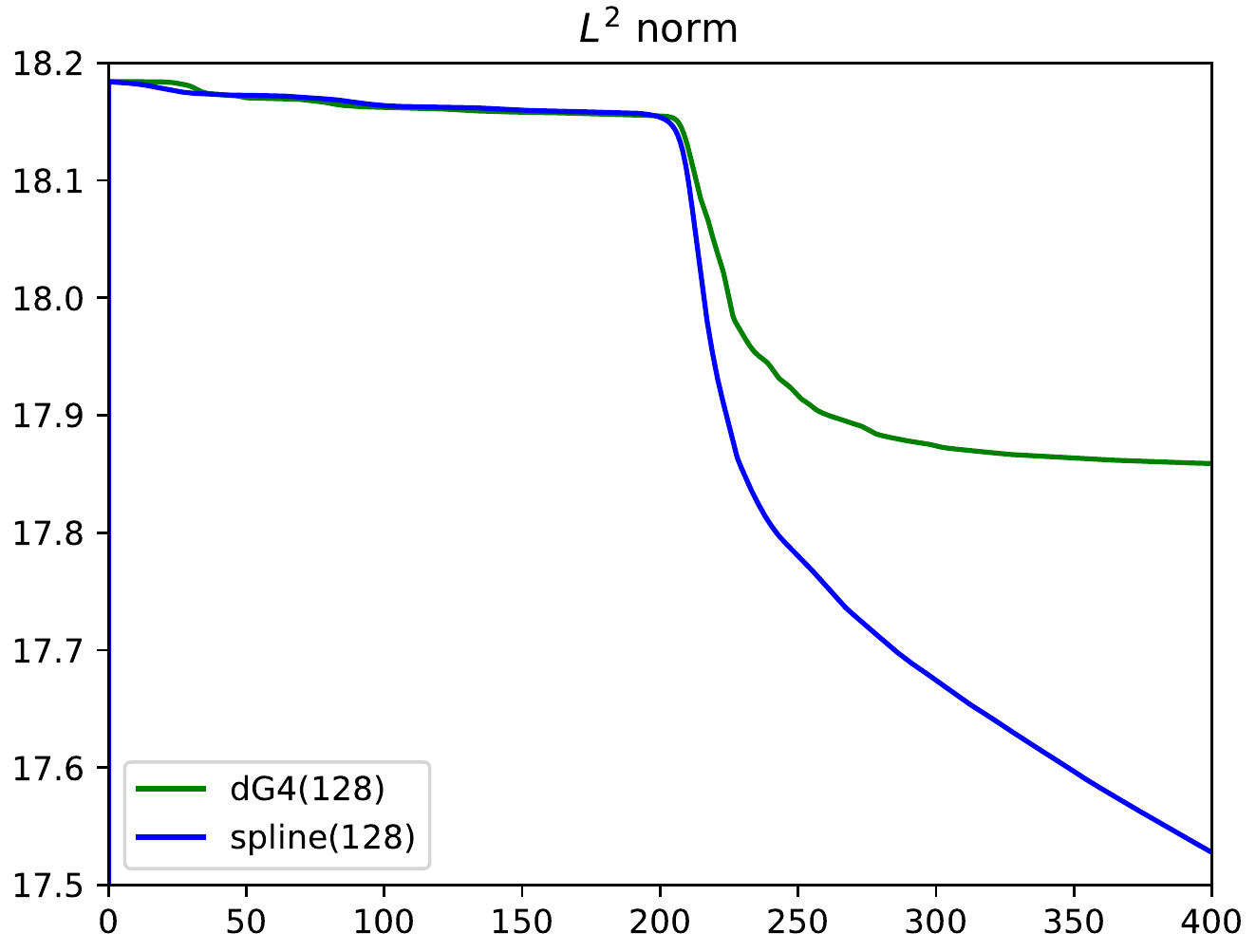}\includegraphics[width=6cm]{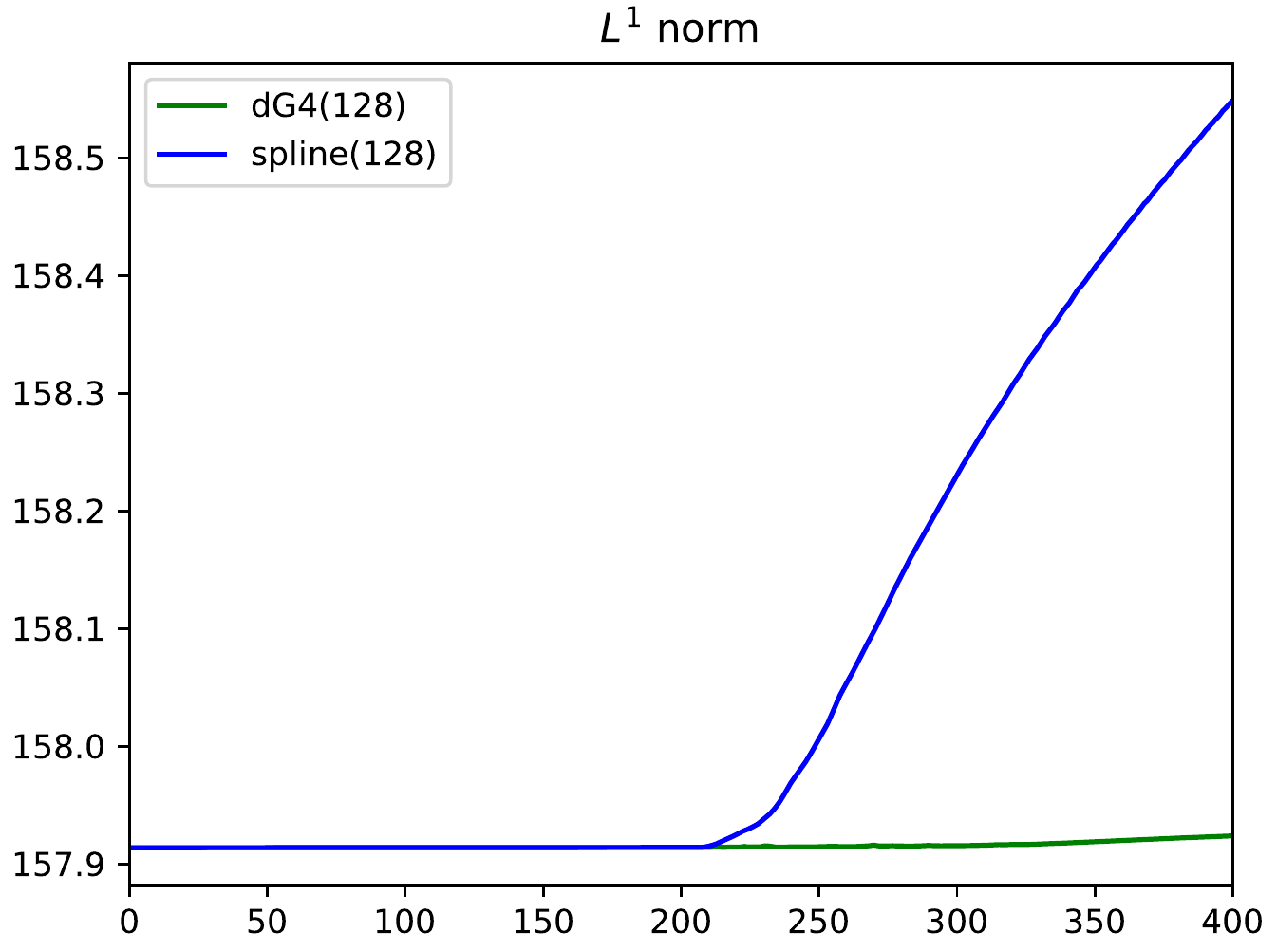}
\par\end{centering}
\caption{The time evolution of the electric energy for the two-stream instability
given by equation (\ref{eq:iv-tsi}) is shown for $32$ (top left),
$64$ (top right), and $128$ (center left) degrees of freedom. In
addition, the energy error (center right), the $L^{2}$ norm (bottom
left), and the $L^{1}$ norm (bottom right) are shown as a function
of time. The reference solution is computed using spline interpolation
and $512$ grid points (per direction).\label{fig:tsi}}
\end{figure}

Figure \ref{fig:tsi} shows a comparison of the fourth order discontinuous
Galerkin method and spline interpolation for $32$, $64$, and $128$
degrees of freedom (per direction). For the dG4 method already $32$
degrees of freedom are sufficient in order to accurately predict the
time and growth of the instability. On the other hand, this is obviously
not true for spline interpolation (which shows an instability that
is too late and does not reach the peak for the electric energy).
While increasing the number of degrees of freedom for spline interpolation
results in significant improvements with respect to these metrics,
the discontinuous Galerkin approach always retains an advantage in
the nonlinear phase. More specifically, spline interpolation exhibits
a decay of the electric energy that is not present in the exact solution.

This behavior can be explained by looking at the time evolution of
the $L^{2}$ norm (see Figure \ref{fig:tsi}). We observe that the
$L^{2}$ norm for dG4 is significantly better conserved compared to
spline interpolation once we enter the nonlinear phase. Thus, more
diffusion is introduced by spline interpolation which results in the
corresponding (unphysical) decrease in electric energy observed (and
the absence of the same artefact for the discontinuous Galerkin approach).

For the present example, the difference between performance in the
asymptotic regime (valid only for small times) and performance for
time scales of physical interest, can be highlighted very clearly.
In fact, if we look at the plots in Figure \ref{fig:tsi}, we see
that in the quiescent phase, the numerical solution obtained with
spline interpolation follows the exact solution somewhat longer. However,
even for $128$ degrees of freedom agreement is only achieved up to
approximately $t=30$ (for spline interpolation). Nevertheless, the
dG4 method is clearly superior in terms of both predicting the onset
of the instability and in providing qualitatively correct predictions
for the nonlinear phase (i.e.~in obtaining the physics of interest).

Let us further remark that the $L^{1}$ norm is significantly better
conserved by the discontinuous Galerkin approach compared to spline
interpolation. Both methods conserve energy quite well (at least up
to $10^{-3}$). Although for the present simulation there is some
advantage for spline interpolation.

We conclude that for the two-stream instability considered here, the
discontinuous Galerkin approach is clearly the superior choice compare
to spline interpolation.

\bibliographystyle{plain}
\bibliography{vlasov-comparison}

\end{document}